\documentclass[a4paper,12pt]{article}
\usepackage{graphicx}
\usepackage{amsmath,amssymb}

\newtheorem{cor}{Corollary}[section]

\newtheorem{theo}{Theorem}[section]
\newtheorem{defi}{Definition}[section]

\newcommand{\be}{\begin{equation}}
\newcommand{\ee}{\end{equation}}

\title{ \bf  NONLINEAR INVARIANTS OF   PLANAR POINT CLOUDS TRANSFORMED BY MATRICES
\rm}

\author{ \parbox{3 in}
{\centering Stelios Kotsios, Evangelos Melas}\\
   \hspace{ 0.5 in}\\
   Faculty of Economics,\\
   Department of Mathematics and Computer Science, \\
National and Kapodistrian University of Athens\\
         Sofokleous 1, Athens 10559, Greece\\
         {\tt \small skotsios@econ.uoa.gr}, {\tt \small emelas@econ.uoa.gr} }

\date{}

\begin{document}
\textwidth 6.50 in
\textheight 9 in

\maketitle

{\bf Abstract:} {\it
The goal of this paper is to present
invariants of planar point clouds, that is functions which take the same value
before and after a linear transformation of a planar point cloud via a $2 \times 2$ invertible matrix.
In the approach we adopt here, 
these invariants are functions of two variables
derived from the  least squares straight line of the planar point cloud under consideration. A linear transformation of
a point cloud induces a nonlinear transformation of
these variables.
The said invariants
are solutions to certain
Partial Differential Equations,
which are obtained by
employing Lie theory. We find cloud invariants 
in the general case of a four$-$parameter transformation 
matrix, as well as, cloud invariants of various one$-$parameter
sets of transformations which can be practically implemented. 
Case studies and simulations which verify our findings are also provided.  }
\section{Introduction}
Analysing level shapes is the key problem in many computer science areas, as image analysis, geometric computing, optical character recognition e.t.c. \cite{kn:FLU,kn:kalouptsidis}.
Usually, by means of the modern sensing technology, we make detailed scans of complex plane objects by generating point cloud data, consisting from thousands or millions of points. Then we study the underlying properties, either by creating appropriate models or by discovering properties which remain constant under
sets of transformations or under the influence of
noise distortions.

\par
\noindent

In particular, when we deal with planar set of points,
a basic approach, which is widely used, is that of determining quantities which can characterize collectively the behaviour of the whole set,
as well as its change, when a transformation is applied to it. In other words, we determine quantities which can represent the planar set of points under consideration, as a whole.
\par \noindent
One approach along this line is the classical work of  Ming$-$Kuei Hu \cite{kn:hu}, who introduced the moment invariants methodology, followed
in the course of time
by many others
\cite{kn:haim, kn:suk, kn:mitra, kn:jan}, to mention but a few. The key element of their approach was
to introduce
the so$-$called moments of planar figures,
in order to identify a planar geometrical figure as a whole,
and then to study their invariants under translation, similitude and orthogonal transformations.
\par \noindent
In
the present
paper, we consider planar set of points,
called henceforth point cloud or cloud of points.
We advocate a different approach, and in order to
characterize collectively
the behaviour of the whole cloud of points
we introduce two variables
$M$ and $H$. These variables stem from
the least squares line assigned to these points. In fact $M$ is the slope of this line, and $H$ is a variation of the $y$$-$intercept of this line.

Any transformation of the cloud of points, by means of a $2 \times 2$ matrix, induces a nonlinear transformation, to be precise a rational one, of the quantities $M$ and $H$. We assume
throughout this paper that any $2 \times 2$ transformation 
matrix of a cloud of points is invertible.
\par \noindent
The purpose
of this paper is to find
cloud invariants, or shortly invariants, expressed in terms
of the variables $M$ and $H$.
By the term ``invariants'' \cite{kn:neusel}
we mean functions which take
the same value at the original and at the transformed values of $M$ and $H$, when a cloud of points undergoes a transformation with a $2 \times 2$ matrix.
In order to solve this problem we use Lie Theory \cite{kn:gil,cant}. As a result, the problem is reduced to solving certain Partial Differential Equations. Any solution to these PDEs, provides us with a cloud invariant.
\par \noindent
Since the problem, in its general form,
does not have a solution which can be practically
implemented, we are examining the special case of an one$-$parameter set of transformations.
This is the case
when the entries of a transformation matrix
are functions of one parameter only.
In this case a general solution is found by using 
Lie theory implemented with symbolic computation.

\par \noindent
At first sight it might seem that restricting
ourselves to a one$-$parameter set of transformations,
useful as it may be, cannot be of great
use. However, this is not the case, because as
we point out in section \ref{linear} \it any \normalfont given matrix belongs to one such
one$-$parameter set of transformations.
As a result we find a family of cloud invariants
for \it any \normalfont given matrix and this certainly lends itself to practical implementation.
\par \noindent
By practical implementation we mean that these invariants
can be used as a tool for studying changes of planar figures
and for creating proper software which monitors and displays these changes in real time. This would have many applications in optical character recognition, as well as, in image analysis and computer graphics techniques; icons
created by the same ``source'' will be readily identified.

\par \noindent
This potential application of our results  suggests a direction for future research.
The cloud of points may come from an icon
which has a parabolic$-$\newline
like shape. In this
case it is natural to look for cloud invariants
which are expressed in terms of variables
which appear as coefficients, or variations
thereof,  of the parabola which is the best
fit for the cloud points we consider.
Comparison with already existing methodologies,
which address the same questions, via simulations and computational experiments, will be also the subject of future research.
\par
\noindent
In section \ref{BQ} we introduce the
variables M and H, and we find the transformation
of these variables which is induced from a transformation with a $2 \times 2$ matrix of the cloud of points under consideration.
In section \ref{invariants1} we define the notion
of an invariant function, and we prove in the
current case, Lie's theory fundamental result that the nonlinear invariant condition is equivalent
to a linear condition provided that the invariant function is properly analytic. Moreover, by using this linear condition, we find the cloud invariants in the general case of a four$-$parameter transformation matrix of the cloud of points under consideration.
In section \ref{onepar2} we find a family of
cloud invariants for a general one$-$parameter
set of transformations. In section \ref{invariants}
we find families of cloud invariants for various sets 
of transformations. We also
find a family  of cloud invariants for
a ``linear'' one$-$parameter set of transformations and we point out that
\it any \normalfont given matrix belongs to such a set.
In section \ref{simulations}  we verify our results with 
simulations and computational experiments in a cloud of 10.000
points. 
In section \ref{remarks} we close the paper with some 
concluding remarks.

\section{The Basic Quantities and their Transformations}
\label{BQ}
In this section we
present
 two quantities $M$ and $H$ which characterize collectively a cloud of points and serve as the independent variables of the invariant functions we are going to construct.
  They originate from the least squares straight line fitted to the cloud of points under consideration.
\par \noindent
Let $(x_i,y_i)$, $i=1,...,N$, be a cloud of points on the plane. We define the quantities:
  \begin{eqnarray}
 	\label{Mdef}
 	M & = & \frac{N\sum_{i=1}^Nx_iy_i-\sum_{i=1}^Nx_i \sum_{i=1}^Ny_i }{N\sum_{i=1}^Nx_i^2-(\sum_{i=1}^Nx_i)^2}, \ \ \rm{and}, \\ 
 	\label{Hdef}
 	H & = & \frac{N\sum_{i=1}^Ny_i^2-(\sum_{i=1}^Ny_i)^2}{N\sum_{i=1}^Nx_i^2-(\sum_{i=1}^Nx_i)^2}.
 \end{eqnarray}
 $M$ is the slope of the least squares straight line and $H$ is suggested by the calculations.
We call them the {\bf linear coefficients} of the cloud.
Sometimes $M$ is referred as the {\bf slope of the cloud} and $H$ as {\bf the constant term} of the cloud.
A transformation of the cloud of points
under the action of a $2 \times 2$ matrix
induces a transformation to $M$ and $H$.
This last transformation  is of prime importance to
our construction of invariant functions and so we proceed
to find it.
Firstly, we need a definition.

\begin{defi}
Let $(x_i,y_i)$, $i=1, \ldots, N$, be a cloud of points and $A=\left(\begin{array}{cc} \alpha & \beta \\ \gamma & \delta \end{array} \right)$,
$\alpha,\beta,\gamma,\delta \in {\bf R}$,  a given $2 \times 2$ matrix. Let us suppose that every point
$(x_i,y_i)$, $i=1, \ldots, N$,
of the cloud undergoes a
transformation $T_{A}: \left(\begin{array}{c}
x_i \\
y_i
\end{array} \right) \to \left(\begin{array}{c}
\hat{x}_i \\
\hat{y}_i
\end{array} \right)$ according to the rule
$\left(\begin{array}{c}
\hat{x}_i \\
\hat{y}_i
\end{array} \right)=A\left(\begin{array}{c}
x_i \\
y_i
\end{array} \right)$.
\par \noindent
We say that the cloud is transformed
under the matrix $A$, and in particular we say that the cloud $(\hat{x}_i,\hat{y}_i)$, $i=1, \ldots, N,$ is the transformation of the cloud $(x_i,y_i)$, $i=1, \ldots, N,$ under the matrix $A$.
\end{defi}
The transformation of $M$ and $H$ induced by
a transformation of the cloud of points via a matrix $A$ is given in the
following Theorem.

\begin{theo}
Let $(x_i,y_i)$, $i=1,2,\ldots,N$, be a cloud of points with linear coefficients $M$ and  $H$.
Let $(\hat{x}_i,\hat{y}_i)$, $i=1,2,\ldots,N$, be the transformation of the cloud $(x_i,y_i)$, $i=1, \ldots, N,$ under
a matrix
$A=\left(\begin{array}{cc} \alpha & \beta \\ \gamma & \delta \end{array} \right)$, $\alpha,\beta,\gamma,\delta \in {\bf R}$.
 Let $\hat{M}$ and  $\hat{H}$ be the linear coefficients of the cloud  $(\hat{x}_i,\hat{y}_i)$, $i=1,2,\ldots,N$.
 Then the following relations hold
\begin{eqnarray}
\label{transfM}
\hat{M}&=&\frac{(\alpha \delta+\beta \gamma)M+\beta\delta H+\alpha \gamma}{2\alpha \beta M+\beta^2 H+\alpha^2}, \\
\label{transfH}
\hat{H}&=&\frac{2\gamma \delta M+\delta^2 H +\gamma^2}{2\alpha \beta M+\beta^2 H+\alpha^2}.
\end{eqnarray}
\end{theo}
{\bf Proof:}
Let $(x_i,y_i)$, $i=1,2,\ldots,N$, be a cloud of points with linear coefficients $M$ and  $H$.
It is  convenient to define the quantities $M_n$, $H_n$, and $D$, as follows
\begin{eqnarray}
M_n & = & N\sum_{i=1}^Nx_iy_i-\sum_{i=1}^Nx_i \sum_{i=1}^Ny_i,  \\
H_n & = & N\sum_iy_i^2-(\sum_iy_i)^2,  \\
D & = &  N\sum_ix_i^2-(\sum_ix_i)^2.
\end{eqnarray}
The relations (\ref{Mdef}) and  (\ref{Hdef}) which define the linear coefficients $M$ and  $H$ of the cloud of points under consideration can now be written in
the following shorter form:
\be
\label{newdef}
M=\frac{M_n}{D}, \quad H=\frac{H_n}{D}.
\ee
A transformation of the cloud of points under a matrix
$A=\left(\begin{array}{cc} \alpha & \beta \\ \gamma & \delta \end{array} \right),$ $\alpha,\beta,\gamma,\delta \in {\bf R},$
induces a transformation to the quantities
$M_n, \ H_n, \ D.$
The induced transformed values $\hat{M}_n, \ \hat{H}_n, \ \hat{D}$, which are assigned to the cloud
$(\hat{x}_i,\hat{y}_i)$, $i=1,2,\ldots,N$, are calculated as follows:
\[ \hat{M}_n=N\sum_{i}\hat{x}_i\hat{y}_i-\sum_{i}\hat{x}_i\sum_{i}\hat{y}_i=\]
\[=N\sum_{i}(\alpha x_i+\beta y_i)(\gamma x_i +\delta y_i)-\sum_{i}(\alpha x_i+\beta y_i)\sum_{i}(\gamma x_i +\delta y_i)=\]
\[=\alpha \gamma \left[ N\sum x_i^2-\left( \sum_{i} x_i \right)^2\right]+(\alpha \delta +\beta \gamma)\left[ N \sum_{i}x_iy_i-\sum_{i}x_i\sum_{i}y_i\right]+\]
\begin{equation}\label{29}
+\beta \delta \left[ N \sum_{i}y_i^2-\left( \sum_{i}y_i^2 \right)\right]=\alpha \gamma D+(\alpha \delta+\beta \gamma)M_n+\beta \delta H_n,
\end{equation}
\[\hat{H}_n=N\sum_{i}\hat{y}_i^2-\left(\sum_{i}\hat{y}_i \right)^2=N\sum_{i}(\gamma x_i+\delta y_i)-\left[\sum_{i}(\gamma x_i+\delta y_i)\right]^2=\]
\begin{equation}\label{30}
=\gamma^2\left[N \sum_{i}x_i^2-\left(\sum_{i}x_i\right)^2\right]+\delta^2\left[ N \sum_{i}y_i^2-\left(\sum_{i}y_i^2\right)\right]+
\end{equation}
\[ + 2\gamma\delta\left[N\sum_{i}x_iy_i-\sum_{i}x_i\sum_{i}y_i \right]=\gamma^2D+\delta^2H_n+2\gamma\delta M_n,\] and,

\[ \hat{D}=N\sum_{i}\hat{x}_i^2-(\hat{x}_i)^2=N\sum_{i}(\alpha x_i +\beta y_i)^2- \left[\sum_{i}(\alpha x_i +\beta y_i) \right]^2= \]
\begin{equation}\label{31}
=\alpha^2 \left[ N \sum_{i}x_i^2-\left(\sum_{i}x_i \right)^2 \right]+\beta^2\left[N \sum_iy_i^2-\left(\sum_{i}y_i^2\right) \right]+
\end{equation}

\[ + 2 \alpha \beta \left[ N \sum_{i}x_iy_i -\sum_{i}x_i \sum_{i}y_i\right]=\alpha^2 D + \beta^2 H_n +2 \alpha \beta M_n. \]
From relations (\ref{newdef}), (\ref{29}),  (\ref{30}), and (\ref{31}), we conclude that  the linear coefficients $\hat{M}$ and  $\hat{H}$ of the cloud  $(\hat{x}_i,\hat{y}_i)$, $i=1,2,\ldots,N$,
are given by the relations (\ref{transfM}), (\ref{transfH}) and the theorem has been proved.
$\square$
\par
\noindent
We denote the set of transformations
(\ref{transfM}) and (\ref{transfH}) by
$\mathcal T(A).$
These are transformations of the form
\be
\hat{M}={\cal M}(M,H,\alpha,\beta,\gamma,\delta),\quad \quad \hat{H}={\cal H}(M,H,\alpha,\beta,\gamma,\delta).
\ee
If $A$ is a matrix with entries ($\alpha,\beta;\gamma,\delta$),
then we can associate to it an element of the set ${\cal T}(A),$
namely the transformation given by (\ref{transfM}) and (\ref{transfH}).
We denote this transformation by ${\cal T}(A)_{(\alpha,\beta,\gamma,\delta)}.$
The following remarks are in order regarding this association:
\begin{itemize}
\item{The set of transformations ${\cal T}(A)$
form a Lie group, under the usual composition of transformations, if and only if the set of matrices  $A$ form also a Lie group, under the usual multiplication of matrices, namely the group $GL(2)$, i.e., the group of $2 \times 2$ invertible matrices. }
\item{This association is not one$-$to$-$one. Indeed, one can easily check that ${\cal T}(A)_{(\alpha,\beta,\gamma,\delta)}={\cal T}(A)_{(\kappa \alpha,\kappa \beta,\kappa \gamma,\kappa \delta)}$,
    $\kappa \in R, \ \kappa \neq 0.$ Therefore \it{all} \normalfont
    matrices  $\kappa A$  are associated to the same element
    ${\cal T}(A)_{(\alpha,\beta,\gamma,\delta)}$
    of ${\cal T}(A).$ }
\item{We can make the association between the sets $A$ and ${\cal T}(A)$ one$-$to$-$one by assigning arbitrarily a fixed non$-$zero value to any of the entries ($\alpha,\beta;\gamma,\delta$) of the matrices of the set $A$.  }
\item{In this paper we prefer not to make this association one$-$to$-$one because this may give the false impression that restrictions are imposed to the set of transformations $A$ which act on the cloud of points under consideration.}

     \item{Needless to say that the results are identical regardless of whether we make or we do not make the association between the sets $A$ and ${\cal T}(A)$ one$-$to$-$one.
    }
\end{itemize}
We note that in the search for invariants we do not need to restrict to the case where $A$, and therefore ${\cal T}(A)$, is a group. As it will become evident from the proof in the next section, and as it will be demonstrated in the example given in subsection \ref{linear}, what it is really necessary is that the set $A$, and subsequently the set ${\cal T}(A)$, must contain the identity element. The identity elements of both   $A$ and ${\cal T}(A)$ are obtained
when $\alpha=1,\beta=0,\gamma=0, \ \rm{and}, \ \delta=1$. For brevity we write
$e=(1,0,0,1)$ and we denote by
$e_i,$ $i=1,2,3,4,$ its components.
\section{Invariants}
\label{invariants1}
A main objective in cloud of points theory is that of finding
invariants.
These are quantities which remain unchanged whenever a cloud of points is transformed under the action of a $2\times 2$ matrix. Invariant quantities enable us
to recognize clouds of points arising from the same ``source''.
\par
\noindent
In our approach the entities which identify a cloud of points are $M$ and $H$. Therefore, we are looking for invariants
which are functions of these two quantities.
This is formalized in the following definition:

\begin{defi}
Let $(x_i,y_i)$, $i=1,...,N,$ be a cloud of points on the plane with linear coefficients $M$ and $H.$
Let $(\hat{x}_i,\hat{y}_i)$, $i=1, \ldots, N,$ be the transformation of the cloud $(x_i,y_i)$, $i=1, \ldots, N,$ under a matrix $A$. Let $\hat{M}$, $\hat{H}$ be the linear coefficients of the cloud $(\hat{x}_i,\hat{y}_i)$, $i=1, \ldots, N.$
$\hat{M}$ and $\hat{H}$ are the transformed values of $M$ and $H$ under the induced
set of
transformations $\mathcal T$$(A).$
We say that a function $I: {\bf R}^2 \to {\bf R}$ is a cloud invariant if and only if

\be
\label{invcon}
I(\hat{M},\hat{H})=I(M,H).
\ee
\noindent
\end{defi}
\noindent
The next theorem is the key result in our study because it 
provides us with a mechanism for finding cloud invariants.
Its proof, is the proof in our case, of Lie's theory fundamental
result \cite{cant} that the nonlinear condition (\ref{invcon}) is equivalent
to a linear condition provided that the invariant function is 
properly analytic. The details are as follows:

\begin{theo}\label{main1}
	Let $(x_i,y_i)$, $i=1,...,N,$ be a cloud of points on the plane with linear coefficients $M$ and $H,$ and let
$(\hat{x}_i,\hat{y}_i)$, $i=1,...,N,$ be the transformation
of the cloud $(x_i,y_i)$, $i=1,...,N,$ under
a matrix $A=\left( \begin{array}{cc}
	\alpha & \beta \\
	\gamma & \delta
	\end{array}
	\right)$. A function

 $I: {\bf R}^2 \to {\bf R}$, analytic in the parameters
  $\alpha, \beta, \gamma,
  \delta$,  is a cloud invariant if and only if the following equations hold simultaneously
\begin{eqnarray}
	\label{con0}
	\xi^{1}_{\alpha} \frac{\partial I}{\partial M} +   \xi^{2}_{\alpha}  \frac{\partial I}{\partial H}&=&0, \\
\label{con1}
\xi^{1}_{\beta} \frac{\partial I}{\partial M} +   \xi^{2}_{\beta}  \frac{\partial I}{\partial H}&=&0, \\
\label{con2}
\xi^{1}_{\gamma} \frac{\partial I}{\partial M} +   \xi^{2}_{\gamma}  \frac{\partial I}{\partial H}&=&0, \\
\label{con3}
\xi^{1}_{\delta} \frac{\partial I}{\partial M} +   \xi^{2}_{\delta}  \frac{\partial I}{\partial H}&=&0,
\end{eqnarray}	
where,
\[ \xi_Q^1=\left(\frac{\partial \hat{M}}{\partial Q}\right )_{e} =\frac{\partial {\cal M}(M,H,1,0,0,1)}{\partial Q},   \xi_Q^2=\left(\frac{\partial \hat{H}}{\partial Q}\right )_{e}  =\frac{\partial {\cal H}(M,H,1,0,0,1)}{\partial Q},\] $Q=\alpha,\beta, \gamma, \delta.$

\end{theo}
\noindent
{\bf Proof:}
Let $(x_i,y_i)$, $i=1,...,N,$ be a cloud of points on the plane with linear coefficients $M$ and $H.$
Let $(\hat{x}_i,\hat{y}_i)$, $i=1, \ldots, N,$ be the transformation of the cloud $(x_i,y_i)$, $i=1, \ldots, N,$ under a matrix $A=\left( \begin{array}{cc}
	\alpha & \beta \\
	\gamma & \delta
	\end{array}
	\right)$. Let $\hat{M}$, $\hat{H}$ be the linear coefficients of the cloud $(\hat{x}_i,\hat{y}_i)$, $i=1, \ldots, N.$
$\hat{M}$ and $\hat{H}$ are the transformed values of $M$ and $H$ under the induced
set of
transformations $\mathcal T$$(A),$
given by equations (\ref{transfM}) and (\ref{transfH}).
Let $I(\hat{M},\hat{H})$  be a real$-$valued function analytic in the
parameters $\alpha, \beta, \gamma, \rm{and} \ \delta.$
The Taylor expansion of $I(\hat{M},\hat{H}),$ with center $e$, reads:
\begin{eqnarray}
\label{T}
I(\hat{M},\hat{H})&=&I(M,H)  +  (\alpha-1)  \left(  \frac{\partial I(\hat{M},\hat{H}) }{\partial \alpha} \right)_{e} +\beta  \left(  \frac{\partial I(\hat{M},\hat{H})}{\partial \beta} \right )_{e} +
\nonumber \\
&& \gamma  \left(  \frac{\partial I(\hat{M},\hat{H})}{\partial \gamma} \right )_{e}+ (\delta-1)  \left(  \frac{\partial I(\hat{M},\hat{H}) }{\partial \delta} \right)_{e} +  \frac{1}{2!} \left ( (\alpha-1)^{2}  \right .
\nonumber \\
&&  \left( \frac{\partial^{2} I(\hat{M},\hat{H})}{\partial \alpha^{2}}\right)_{e}  +  \beta^{2} \left( \frac{\partial^{2} I(\hat{M},\hat{H})}{\partial \beta^{2}} \right )_{e}+ \gamma^{2}\left( \frac{\partial^{2} I(\hat{M},\hat{H})}{\partial \gamma^{2}} \right)_{e} + \nonumber \\
&&(\delta-1)^{2} \left( \frac{\partial^{2} I(\hat{M},\hat{H})}{\partial \delta^{2}} \right )_{e}+2 (\alpha-1)\beta  \left(  \frac{\partial I(\hat{M},\hat{H})}{\partial \alpha} \right )_{e}
+ \cdots \nonumber \\
&& \left .
2\gamma (\delta-1) \left(  \frac{\partial I(\hat{M},\hat{H})}{\partial \gamma} \right )_{e}
\left(  \frac{\partial I(\hat{M},\hat{H})}{\partial \delta} \right )_{e} \right) +\cdots \ .
\end{eqnarray}

\noindent
The form of the functional dependence of $I(\hat{M},\hat{H})$ on the parameters $\alpha,\beta,\gamma,\delta$ allows to simplify (\ref{T})
in a way which suggests the conclusion of the theorem.
To illustrate  the point at hand we use the chain rule to evaluate the derivative $\displaystyle \left ( \frac{\partial I(\hat{M},\hat{H})}{\partial \alpha} \right )_{e}:$

\be
\left(  \frac{\partial I(\hat{M},\hat{H})}{\partial \alpha} \right )_{e}
=
\left(  \frac{\partial I(\hat{M},\hat{H})}{\partial \hat{M}} \frac{\partial \hat{M}}{\partial \alpha} +\frac{\partial I(\hat{M},\hat{H})}{\partial \hat{H}} \frac{\partial \hat{H}}{\partial \alpha}\right )_{e}.
\ee
\be
\label{EQ}
\hspace{-0.4cm}
\left(  \frac{\partial I(\hat{M},\hat{H})}{\partial \alpha} \right )_{e}
 =
 \left(  \frac{\partial I(\hat{M},\hat{H})}{\partial \hat{M}}\right )_{e}  \left(\frac{\partial \hat{M}}{\partial \alpha}\right )_{e} + \left(\frac{\partial I(\hat{M},\hat{H})}{\partial \hat{H}} \right )_{e} \left(\frac{\partial \hat{H}}{\partial \alpha}\right )_{e}.
 \ee
By introducing the quantities:
\be
\xi_Q^1=\left(\frac{\partial \hat{M}}{\partial Q}\right )_{e}, \
\xi_Q^2=\left(\frac{\partial \hat{H}}{\partial Q}\right )_{e}, \
Q=\alpha,\beta, \gamma, \delta,
\ee
and by noting
\be
\left(  \frac{\partial I(\hat{M},\hat{H})}{\partial \hat{M}}\right )_{e}= \frac{\partial I(M,H)}{\partial {M}}, \  \rm{and}, \normalfont \
\left(  \frac{\partial \it{I}(\hat{M},\hat{H})}{\partial \hat{H}}\right )_{\it{e}}= \frac{\partial \it{I}(M,H)}{\partial \it{H}},
\ee
we can rewrite equation (\ref{EQ}) in the following shorter form
\be
\label{EQ1}
\left( \frac{\partial I(\hat{M},\hat{H})}{\partial \alpha}\right )_{e} =   \xi^1_\alpha\frac{\partial I(M,H)}{\partial {M}}  + \xi^2_\alpha\frac{\partial I(M,H)}{\partial {H}}.
\ee

By introducing the
operator
\be
X_{\alpha} = \xi^{1}_{\alpha}  \frac{\partial}{\partial M}+\xi^{2}_{\alpha}  \frac{\partial}{\partial H},
\ee
we rewrite equation (\ref{EQ1}) as
\be
\left ( \frac{\partial I(\hat{M},\hat{H})}{\partial \alpha}
\right )_{e}
=X_\alpha I,
\ee
where for short we wrote $I$ instead of $I(M,H).$
For the second order derivative $\left ( \frac{\partial^{2} I(\hat{M},\hat{H})}{\partial \alpha^{2}} \right )_{e}$ we have:
\begin{eqnarray}
\label{sdn}
 \left ( \frac{\partial^{2} I(\hat{M},\hat{H})}{\partial \alpha^{2}} \right )_{e} & = &
 \left ( \frac{\partial \left ( \frac { \partial I(\hat{M},\hat{H})} {\partial \alpha } \right )
 }{\partial \alpha} \right )_{e}
 =
\xi^1_\alpha\frac{\partial \left(\frac{\partial I(\hat{M},\hat{H})}{\partial \alpha}\right)_{e}}{\partial {M}}  + \xi^2_\alpha\frac{\partial \left(\frac{\partial I(\hat{M},\hat{H})}{\partial \alpha}\right)_{e}}{\partial {H}}
\nonumber \\
& = & X_{\alpha}   \left( X_{\alpha} I  \right).
\end{eqnarray}
A similar analysis applies to the derivatives of all orders in the Taylor expansion (\ref{T}).
As a result the Taylor expansion (\ref{T}) reads:
\begin{eqnarray}
\label{Lie8}
\hspace{-1cm}
I(\hat{M},\hat{H})   &=&   I(M,H) + \sum_{i=1}^{4}(\mathcal{Q}_{i}-e_i)  \left( X_{\mathcal Q_i} I   \right)
+ \nonumber \\
&& \hspace{-1cm} \frac{1}{ 2!} \sum_{i,j=1}^{4}(\mathcal{Q}_{i}-e_i)(\mathcal{Q}_{j}-e_j)
 X_{\mathcal Q_i} \left( X_{\mathcal Q_j} I  \right) + \nonumber \\
&& \hspace{-1cm} \frac{1}{ 3!} \sum_{i,j,k=1}^{4}(\mathcal{Q}_{i}-e_i)(\mathcal{Q}_{j}-e_j)
(\mathcal{Q}_{k}-e_k) X_{\mathcal Q_i} \left( X_{\mathcal Q_j}  \left( X_{\mathcal Q_k} I  \right)  \right)
 + \cdots \ .
\end{eqnarray}
For convenience, by $\mathcal{Q}$ we denote the vector $(\alpha,\beta, \gamma, \delta),$ and by  $\mathcal{Q}_i,$ $i=1,2,3,4,$ its components.

\noindent
From equation (\ref{Lie8})  we conclude that when the \it  linear infinitesimal \normalfont conditions
\be \label{lin}  X_{\mathcal Q_i}I =0, \quad \mathcal Q_i=\alpha,\beta,\gamma,\delta,
\ee
are satisfied then $I(\hat{M},\hat{H})=I(M,H)$.
Therefore $I$ is a cloud invariant.
\noindent
Conversely, when $I$ is a cloud invariant,
then $I(\hat{M},\hat{H})=I(M,H),$
and equation (\ref{Lie8}) gives:

 \[  \hspace{-2cm}   \sum_{i=1}^{4}(\mathcal{Q}_{i}-e_i)  \left( X_{\mathcal Q_i} I   \right)
+ \frac{1}{ 2!} \sum_{i,j=1}^{4}(\mathcal{Q}_{i}-e_i)(\mathcal{Q}_{j}-e_j)
X_{\mathcal Q_i} \left( X_{\mathcal Q_j} I  \right) +\] \begin{eqnarray}
\label{Lie9}
\hspace{-2cm}
&& \frac{1}{ 3!} \sum_{i,j,k=1}^{4}(\mathcal{Q}_{i}-e_i)(\mathcal{Q}_{j}-e_j)
(\mathcal{Q}_{k}-e_k) X_{\mathcal Q_i} \left( X_{\mathcal Q_j}  \left( X_{\mathcal Q_k} I  \right)  \right)
+ \cdots =0.
\end{eqnarray}
For every
pair of values
$M$ and $H$
equation (\ref{Lie9}) becomes a polynomial in the variables $\alpha,\beta,\gamma,\delta$. Consequently
equation (\ref{Lie9}) can only hold if for every pair of values
$M$ and $H$
the coefficients of the polynomial vanish, i.e., if the
following relations hold
\begin{eqnarray}
\label{Lie10}
X_{\mathcal Q_i} I & = &
X_{\mathcal Q_i} \left( X_{\mathcal Q_j} I \right ) = \nonumber \\
&&X_{\mathcal Q_i} \left( X_{\mathcal Q_j}   \left( X_{\mathcal Q_k} I \right ) \right )=\cdots =0, \quad i,j,k \in\{1,2,3,4\},
\end{eqnarray}
for every pair
of values
$M$ and $H.$
If $
X_{\mathcal Q_i}I =0, \quad \mathcal Q_i=\alpha, \beta,\gamma,\delta,
$ the rest of the relations (\ref{Lie10}) follow.
Equations  (\ref{lin}),
$  X_{\mathcal Q_i}I =0, \quad \mathcal Q_i=\alpha,\beta,\gamma,\delta,
$
are nothing but equations
(\ref{con0}), (\ref{con1}), (\ref{con2}), and
(\ref{con3}), respectively.

\noindent
This completes the proof. $\square$

\noindent
Sophus Lie's great advance was to replace the
complicated, nonlinear finite invariance condition (\ref{invcon}) by the vastly more useful linear
infinitesimal condition (\ref{lin}) and to recognize that if a function satisfies the infinitesimal condition then it also satisfies
the finite condition, and vice versa, provided that the function is analytic in the parameters
$\alpha,\beta,\gamma,$ and  $\delta.$

\noindent
It is to be noted
that in the proof of
Lie's main Theorem (\ref{main1})
we used the following:
\begin{enumerate}

\item{The assumption that  cloud invariant $I$ is a function analytic in the parameters $\alpha,\beta,\gamma,$ and  $\delta.$}

\item{The assumption that the set of transformations $\mathcal T(A)$, and
    subsequently
    the set of transformations $A$, contain the identity
    element, which is  obtained when $\alpha=1, \ \beta=0,
    \ \gamma=0,$ and $\delta=1.$}

\item{The chain rule for the differentiation of composite functions.}

\end{enumerate}
Nowhere in the proof of
Lie's main Theorem (\ref{main1})
is the assumption made that the
set of transformations $\mathcal T(A)$
is closed under the usual composition of
transformations, or equivalently,
that the
associated set of matrices $A$ is closed
under the usual matrix multiplication.
This will become evident and
exemplified  in   subsection \ref{linear}
where we find cloud invariants $I$
under a
set of transformations $\mathcal T(A)$
which are such that the associated
set of matrices $A$ are not closed
under the usual  multiplication of matrices.

\subsection{Cloud invariants in the general case}

\par
\noindent
As a first application of   Theorem  (\ref{main1})  we find the cloud invariants under a general matrix $A.$ This is the content of the next
Corollary.
\begin{cor}
\label{general}
	If a cloud of points is transformed via a
	matrix $A=\left( \begin{array}{cc}
	\alpha & \beta \\
	\gamma & \delta
	\end{array}
	\right)$, then the only cloud invariants are:
\begin{enumerate}
\item{The constant functions $I(M,H)=c$, $c \in R.$}
\item{The level curve $I(M,H)=0$  of the function $\displaystyle I(M,H)=\frac{H}{M^{2}} - 1$.}
\end{enumerate}
\end{cor}
{\bf Proof:}
\noindent
According to Theorem (\ref{main1}) a function
$I(\hat{M},\hat{H}),$ analytic in the parameters
$\alpha,\beta,\gamma,\delta,$ is a cloud
invariant if and only if it satisfies
the  system  of PDEs   (\ref{con0}), (\ref{con1}), (\ref{con2}), and (\ref{con3}), which read:

\begin{eqnarray}
\xi^{1}_{\alpha}(\bold x)  \frac{\partial I }{\partial M} + \xi^{2}_{\alpha}(\bold x)  \frac{\partial I}{\partial H} & = &
-  M \frac{\partial I}{\partial M} - 2 H  \frac{\partial I}{\partial H}=0, \\
\xi^{1}_{\beta}(\bold x)  \frac{\partial I }{\partial M} + \xi^{2}_{\beta}(\bold x)  \frac{\partial I}{\partial H} & = &
(H - 2 M^{2}) \frac{\partial I}{\partial M} - 2 H M  \frac{\partial I}{\partial H}=0, \\
\xi^{1}_{\gamma}(\bold x)  \frac{\partial I}{\partial M} + \xi^{2}_{\gamma}(\bold x)  \frac{\partial I}{\partial H} & = &
\frac{\partial I}{\partial M} - 2  M  \frac{\partial I}{\partial H}=0, \\
\xi^{1}_{\delta}(\bold x)  \frac{\partial I }{\partial M} + \xi^{2}_{\delta}(\bold x)  \frac{\partial I}{\partial H}& = &
M \frac{\partial I}{\partial M} + 2  H  \frac{\partial I}{\partial H}=0.
\end{eqnarray}

We easily find that the only solutions to the last system of equations are:
\begin{enumerate}
\item{The constant functions $I(M,H)=c$, $c \in R.$}
\item{The level curve $I(M,H)=0$  of the function $\displaystyle I(M,H)=\frac{H}{M^{2}} - 1$.}
\end{enumerate}
This completes the proof.$\square$

\noindent
The second invariant
implies in particular that when the  values of M and H are such that $H=M^{2}$,
then their  transformed values $\hat{H}$ and $\hat{M}$
are
such that $\hat{H}=\hat{M}^{2}$.

\section{The One$-$Parameter Case}
\label{onepar2}

The cloud invariants under a general matrix
$A,$ given in Corollary \ref{general}, do not
lend themselves to practical implementation.
This leads us to examining
particular cases of $A.$
We start by considering the one$-$parameter case,
i.e. the case where the entries of the matrix $A$ are analytic functions of a single parameter $\varphi$. Interestingly enough it turns out that
in this case we can find cloud invariants
in closed form which can be practically implemented. The first step to prove this
assertion is the next Corollary.

\begin{cor}\label{main}
Let $(x_i,y_i)$, $i=1,...,N,$ be a cloud of points on the plane with linear coefficients $M$ and $H,$ and let
$(\hat{x}_i,\hat{y}_i)$, $i=1,...,N,$ be the transformation
of the cloud $(x_i,y_i)$, $i=1,...,N,$ under
a matrix
$A(\varphi)=\left( \begin{array}{cc}
	\alpha(\varphi) & \beta(\varphi) \\
	\gamma(\varphi) & \delta(\varphi)
	\end{array}
	\right)$, where $\alpha(\varphi)$, $\beta(\varphi)$,
 $\gamma(\varphi)$, and  $\delta(\varphi),$ are real analytic functions of a parameter $\varphi \in {\bf R}$. We assume that
 there exists a value of $\varphi$, denoted by $\varphi^*$, such that $A(\varphi^*)={\bf I}$, ${\bf I}$ the $2\times 2$ identity matrix. An analytic function $I: {\bf R}^2 \to {\bf R}$ is a cloud invariant if and only if the next equation holds:
\be \label{onepar1}
[(H-2M^2)\beta'(\varphi^*)-
\delta
M+\gamma'(\varphi^*)]\frac{\partial I}{\partial M}+2[\gamma'(\varphi^*) M-\beta'(\varphi^*) HM- \delta H]\frac{\partial I}{\partial H}=0,
\ee
where, $\delta=\alpha'(\varphi^*)-\delta'(\varphi^*).$

\end{cor}
{\bf Proof:}
In order to find the cloud invariants
we apply Theorem \ref{main1}. The key point
is that in the case under consideration
\it all \normalfont the entries of the
matrix $A(\varphi)$  are functions of
a single parameter $\varphi.$ This
implies in particular that
equations (\ref{con0}), (\ref{con1}),
(\ref{con2}), and (\ref{con3}),
whose solution space are the
cloud invariants, reduce to one equation:
\be
\label{redeq}
\xi^{1}_{\varphi} \frac{\partial I}{\partial M} +   \xi^{2}_{\varphi}  \frac{\partial I}{\partial H}=0,
\ee
$$
\xi_\varphi^1=\left(\frac{\partial \hat{M}}{\partial \varphi}\right )_{e} =\frac{\partial {\cal M}(M,H,1,0,0,1)}{\partial \varphi},   \xi_\varphi^2=\left(\frac{\partial \hat{H}}{\partial \varphi}\right )_{e}  =\frac{\partial {\cal H}(M,H,1,0,0,1)}{\partial \varphi}.$$
Differentiation is now with respect to
$\varphi,$ that is $Q=\varphi.$
Consequently cloud invariants
in the one$-$parameter case are solutions
to equation (\ref{redeq}) which reads:
$$
	[(H-2M^2)\beta'(\varphi^*)-
\delta
M+\gamma'(\varphi^*)]\frac{\partial I}{\partial M}+2[\gamma'(\varphi^*) M-\beta'(\varphi^*) HM- \delta H]\frac{\partial I}{\partial H}=0,$$
where, $\delta=\alpha'(\varphi^*)-\delta'(\varphi^*).$
This completes the proof. $\square$

\par
\noindent
It is difficult to obtain in closed form the whole set of solutions of equation (\ref{onepar1}).
However we can find, in closed form, a wide
subclass of solutions of equation (\ref{onepar1}).
This is the content of the following Theorem.
\begin{theo}\label{stathero}
	A class of solutions of  equation (\ref{onepar1}), and hence a family of invariants
of a
cloud of points $(x_i,y_i)$, $i=1,...,N,$
when it is transformed
under a matrix
$A(\varphi)=\left( \begin{array}{cc}
	\alpha(\varphi) & \beta(\varphi) \\
	\gamma(\varphi) & \delta(\varphi)
	\end{array}
	\right)$,
is given by:
	\begin{equation}\label{coninv}
	I(M,H)=F\left(\frac{M^2-H}{(H
\beta'(\varphi^*) -\gamma'(\varphi^*) + \delta M)^2}\right),
	\end{equation}	
where $F(.)$, is an arbitrary real valued function and $\delta=\alpha'(\varphi^*)-\delta'(\varphi^*)$.	We assume that $\alpha(\varphi)$, $\beta(\varphi)$,
 $\gamma(\varphi)$, and  $\delta(\varphi),$
are real analytic functions of a single parameter $\varphi \in {\bf R}$. We also assume that
 there exists a value of $\varphi$, denoted by $\varphi^*$, such that $A(\varphi^*)={\bf I}$, ${\bf I}$ being the $2\times 2$ identity matrix.	


\end{theo}
{\bf Proof:} In order to find solutions of  equation (\ref{onepar1}) we use the undetermined coefficients method.
This method consists in seeking for solutions $I(M,H)$ of the form:
\be
\label{eq}
F\left(\frac{\sum_{i=1}^{n}\sum_{j=1}^{m}a_{ij}M^iH^j}{\sum_{i=1}^{n}\sum_{j=1}^{m}b_{ij}M^iH^j}\right),
\ee
where $a_{ij}$ and $b_{ij}$ are unknown
coefficients to be determined. By substituting
this particular form of the solution into
equation (\ref{onepar1}) we obtain that a
polynomial in the two variables $M$ and $H$
is equal to zero. The resulting condition,
the coefficients of the polynomial are equal to
zero, gives solution (\ref{coninv}). This completes the proof.
$\square$

\par \noindent
A Corollary
of the previous theorem
is that we can find cloud invariants, when the cloud is transformed under a matrix, provided the matrix is an element of the one$-$parameter set
of transformations $A(\varphi)$ considered in this Theorem.
\begin{cor}\label{arbitrary}
Let $(x_i,y_i)$, $i=1,...,N,$ be a cloud of points and let this cloud be transformed under
a matrix $A=\left(\begin{array}{cc}
	a_{11} & a_{12} \\
	a_{21} & a_{22}
	\end{array}\right) $, $a{ij} \in {\bf R}$.
If there exist real valued, analytic, functions
$\alpha(\varphi),\beta(\varphi),\gamma(\varphi),\delta(\varphi)$ and values $\varphi^*, \varphi_1$, such that:
\begin{enumerate}
		
		\item $\alpha(\varphi^*)=1,\beta(\varphi^*)=0,\gamma(\varphi^*)=0,\delta(\varphi^*)=1$
		
		
		\item 	$\alpha(\varphi_1)=a_{11},\beta(\varphi_1)=a_{12},\gamma(\varphi_1)=a_{21},\delta(\varphi_1)=a_{22}$
	\end{enumerate}
then, the quantity
$$
	I(M,H)=F\left(\frac{M^2-H}{(H
\beta'(\varphi^*) -\gamma'(\varphi^*) + \delta M)^2}\right),
$$
where $F(.)$ is an arbitrary real valued function and $\delta=\alpha'(\varphi^*)-\delta'(\varphi^*),$
	is a cloud invariant. 	
\end{cor}
{\bf Proof:}
This is an immediate consequence of Theorem
\ref{stathero}. $\square$



\section{Examples of Invariants}
\label{invariants}
In this section, by using Theorem theorem \ref{stathero}, we find cloud invariants when a cloud
is transformed under various sets of transformations
$A(\varphi).$ As we pointed out in section \ref{BQ} it is not
necessary the set $A(\varphi)$ to form a group under
the usual multiplication of matrices. Firstly we
consider sets of transformations
$A(\varphi)$ which do form a group and then
we consider a set
$A(\varphi)$ which does not form a group.
Finally, in the last subsection, by using the
previous findings, we find cloud invariants
for any given matrix.

\subsection{Sets $A(\varphi)$ which form a group}
We  start  with simpler sets of transformations $A(\varphi)$ and   gradually proceed to more general cases.

\subsubsection{$A(\varphi)$ is a diagonal matrix}
\label{sub1}
We start by assuming that  $A(\varphi)$ is diagonal and has  the form:
\be
A(\varphi)= \left(\begin{array}{cc}1 & 0 \\ 0 & \varphi \end{array} \right),
\ee
\noindent
where $\varphi \in {\bf R}$. In this case we can easily see that $\varphi^*=1$ and that
$\beta'(\varphi^*)=0, \ \gamma'(\varphi^*)=0,  \ \delta=\alpha'(\varphi^*)-\delta'(\varphi^*)=-1.$
Therefore, according to Theorem \ref{stathero}, a family of cloud invariants is:
\be
F\left(\frac{M^2-H}{(- M)^2}\right)=h\left(\frac{H}{M^2}\right),
\ee
where $F(\cdot)$ and  $h(\cdot)$ are arbitrary real valued functions.

\subsubsection{$A(\varphi)$ is an upper triangular matrix}
\label{sub2}
We assume  that  $A(\varphi)$ is upper triangular and has the form:
\be
A(\varphi)=\left(\begin{array}{cc}1 & \varphi \\ 0 & 1 \end{array}  \right),
\ee
$\varphi \in {\bf R}$.
\noindent
In this case we  easily verify that $\varphi^*=0$ and that
$\beta'(\varphi^*)=1, \ \gamma'(\varphi^*)=0,  \ \delta=\alpha'(\varphi^*)-\delta'(\varphi^*)=0.$
Consequently, according to Theorem \ref{stathero},
a family of cloud invariants is:
\be
F\left(\frac{M^2-H}{(H\cdot 1)^2}\right)=F\left(\frac{M^2-H}{H^2}\right),
\ee
where $F(\cdot)$ is an arbitrary real valued function.

\subsubsection{$A(\varphi)$ is a lower triangular matrix}
\label{sub3}
We assume  that  $A(\varphi)$ is lower triangular and has the form:
\be
A(\varphi)=\left(\begin{array}{cc}1 & 0 \\ \varphi & 1 \end{array}  \right),
\ee
$\varphi \in {\bf R}$.
In this case we  easily find that $\varphi^*=0$ and that
$\beta'(\varphi^*)=0, \ \gamma'(\varphi^*)=1,  \ \delta=\alpha'(\varphi^*)-\delta'(\varphi^*)=0.$
Consequently, according to Theorem \ref{stathero},
a family of cloud invariants is:
\be
F\left(\frac{M^2-H}{(-1)^2}\right)=h(H-M^2),
\ee
where $F(\cdot)$ and  $h(\cdot)$ are arbitrary real valued functions.
\noindent

\subsubsection{$A(\varphi)$ is a rotation matrix}
\label{sub4}
Finally, we assume  that  $A(\varphi)$ is a rotation matrix and has the form:
\be
A(\varphi)=\left(\begin{array}{cc}
\cos \varphi & \sin \varphi \\
-\sin \varphi & \cos \varphi
\end{array} \right),
\ee
$\varphi \in {\bf R}$.
In this case we  easily obtain that $\varphi^*=0$ and that
$\beta'(\varphi^*)=1, \ \gamma'(\varphi^*)=-1,  \ \delta=\alpha'(\varphi^*)-\delta'(\varphi^*)=0.$
Consequently, according to Theorem \ref{stathero},
a family of cloud invariants is:
\be
 F\left(\frac{M^2-H}{(H +1)^2}\right),
 \ee
where $F(\cdot)$ is an arbitrary real valued function.

\subsection{A set $A(\varphi)$ which does not form a group}

\subsubsection{A ``linear" matrix}
\label{linear}
A set of transformations $A(\varphi)$ which subsumes the
sets of transformations considered in subsections
\ref{sub1}, \ref{sub2}, and \ref{sub3} is the set
of ``linear" matrices
\small
\begin{equation}\label{lineartr}
A(\varphi)=\left(\begin{array}{cc}
a_{11} & a_{12} \\
a_{21} & a_{22}
\end{array} \right)+\varphi \left(\begin{array}{cc}
b_{11} & b_{12} \\
b_{21} & b_{22}
\end{array}\right)
=
\left(\begin{array}{cc}
a_{11}+b_{11} \varphi & a_{12}+b_{12} \varphi \\
a_{21}+b_{21} \varphi & a_{22}+b_{22} \varphi
\end{array} \right),
\end{equation}
\normalsize
where $a_{ij},b_{ij} \in {\bf R}$,
and $\varphi$ is a real free parameter.
The set of matrices $A(\varphi^*)$ does not, in general, form a group
under matrix multiplication.
However, we assume that
there exists a value of $\varphi$, denoted by $\varphi^*$, such that $A(\varphi^*)={\bf I}$, ${\bf I}$ being the $2\times 2$ identity matrix.
One can easily check that such a value $\varphi^*$ exists
if and only if the entries $a_{ij},b_{ij}$ satisfy one of the following conditions:
\be
b_{22}\neq 0\land a_{21}=\frac{\left(a_{22}-1\right) b_{21}}{b_{22}}\land a_{12}=\frac{\left(a_{22}-1\right) b_{12}}{b_{22}}\land a_{11}=\frac{a_{22} b_{11}-b_{11}+b_{22}}{b_{22}},
\ee
or
\be
b_{22}=0\land a_{22}=1\land b_{21}\neq 0\land a_{12}=\frac{a_{21} b_{12}}{b_{21}}\land a_{11}=\frac{a_{21} b_{11}+b_{21}}{b_{21}},
\ee
or
\be
b_{22}=0\land b_{21}=0\land a_{22}=1\land a_{21}=0\land b_{12}\neq 0\land a_{11}=\frac{a_{12} b_{11}+b_{12}}{b_{12}},
\ee
or
\be
b_{22}=0\land b_{21}=0\land b_{12}=0\land a_{22}=1\land a_{21}=0\land a_{12}=0\land b_{11}\neq 0.
\ee
In this case we easily find that $\beta'(\varphi^*)=b_{12}, \ \gamma'(\varphi^*)=b_{21},  \ \delta=\alpha'(\varphi^*)-\delta'(\varphi^*)=b_{11}-b_{22}.$
 According to Theorem \ref{stathero}, a family of cloud invariants is:

\begin{equation}\label{grammiko}
F\left(\frac{M^2-H}{(Hb_{12} -b_{21} +(b_{11}-b_{22}) M)^2}\right),
\end{equation}
where $F(\cdot)$ is an arbitrary real valued function.

\subsubsection{Cloud invariants for an arbitrary matrix}
For \it {any} \normalfont
 given matrix there always exists a set of
transformations $A(\varphi),$ of the form  (\ref{lineartr}),
which contains this matrix. In fact one can easily prove
that there exists a two parameter family of such sets
$A(\varphi).$
According to Corollary \ref{arbitrary}, a family of cloud invariants, when a cloud of
points is transformed under this matrix, is given by relation
(\ref{grammiko}).
As a case study we consider the matrix: \be
\mathcal M= \left(
\begin{array}{cc}
0.4 & -0.4 \\
-0.05 & 0.9 \\
\end{array}
\right). \ee
Let $A(\varphi)$ be the set of transformations:
\be
A(\varphi)=\left(\begin{array}{cc}
-2 & -2 \\
-1/4 & 1/2
\end{array} \right)+\varphi \left(\begin{array}{cc}
12 & 8 \\
1 & 2
\end{array} \right)=\left(
\begin{array}{cc}
-2+12 \varphi  & -2+8\varphi \\
-1/4+ \varphi & 1/2+2 \varphi  \\
\end{array}
\right).
\ee
One can easily check that $A(0.2)=\mathcal M.$
We have $b_{12}=8, \ b_{21}=1, \ b_{11}-b_{22}=10.$
 According to Corollary \ref{arbitrary}, a family of cloud invariants is:
	\begin{equation}\label{casestudy}
	F\left(\frac{M^2-H}{(8H -1 + 10 M)^2}\right),
	\end{equation}	
where $F(\cdot)$ is an arbitrary real valued function.
Since there exists a two parameter family of  sets
$A(\varphi),$ of the form (\ref{lineartr}), which contain
$\mathcal M,$ there exists a two parameter family of cloud
invariants of the form (\ref{casestudy}), when a cloud
is transformed via $\mathcal M.$ However,
the explicit form of this two parameter family of cloud
invariants
is not needed
here.

\section{Simulations}
\label{simulations}

\begin{figure}
	\centering
	\includegraphics[scale=0.4]{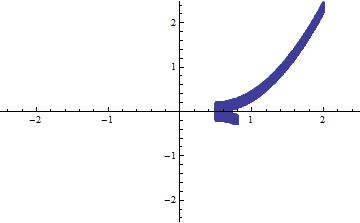}
	\caption{The Original Scheme}
	\label{sxhma0}
\end{figure}

To see how the above theory works in practise, we
 consider a cloud of $10000$ points forming the scheme
of  Figure \ref{sxhma0}. Using  relations (\ref{Mdef})
 and (\ref{Hdef}), we calculate
 the linear coefficients  $M$ and $H$ of the cloud. We find
$M=1.52244  $ and $H=2.46998$.
We transform now this cloud
by using various matrices.
\par
\begin{figure}
	\centering
	\includegraphics[scale=0.4]{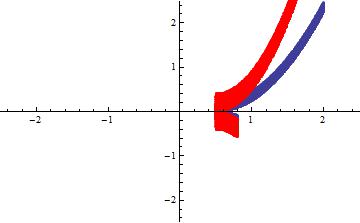}
	\caption{A Diagonal Transformation}
	\label{sxhmadiag}
\end{figure}

 \noindent
 Firstly we let the
the diagonal matrix
$A=\left(\begin{array}{cc}
1 & 0 \\
0 & 2
\end{array} \right)$ to act on the cloud.
Both the initial and the transformed schemes
are depicted in Figure \ref{sxhmadiag}.
Using relations (\ref{transfM}) and (\ref{transfH})
we calculate the linear coefficients of the new cloud and we  obtain $\hat{M}=3.04488$ and $\hat{H}=9.87991$.
According to our findings in subsection \ref{sub1},
a family of cloud invariants, when a cloud is
transformed with the diagonal matrix $A,$ is $I(M,H)= F(H/M^2),$ where $F(\cdot)$ is an
arbitrary real valued function.
Indeed, for the initial and
the transformed linear coefficients,  we find $H/M^2=\hat{H}/\hat{M}^{2}=1.06564.$ It follows
that we have $I(M,H)= I(\hat{M},\hat{H})  $ for \it any \normalfont real valued function $F(\cdot)$.

\par
\begin{figure}
	\centering
	\includegraphics[scale=0.4]{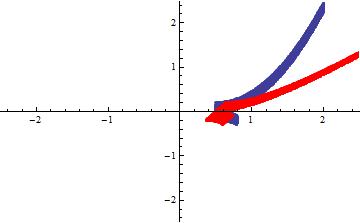}
	\caption{An Upper Triangular Transformation}
	\label{sxhmaupper}
\end{figure}

\noindent
As a second example of transformation,
we let the upper triangular matrix $B=\left(\begin{array}{cc}
1 & 0.7 \\
0 & 1
\end{array} \right)$ to act on the cloud.
The result of this transformation is given
in Figure \ref{sxhmaupper}. The transformation
of the cloud we consider under $B$ has linear coefficients
$\hat{M}=0.748882$ and $\hat{H}=0.568896$.
As we found in subsection \ref{sub2},
a family of cloud invariants, when a cloud is
transformed with the upper triangular matrix $B,$ is $I(M,H)= F((M^2-H)/H^2),$ where $F(\cdot)$ is an arbitrary real valued function. Indeed, we have
$(M^2-H)/H^2 =(\hat{M}^2-\hat{H})/\hat{H}^2 =-0.0249396.$ Consequently
 we have $I(M,H)= I(\hat{M},\hat{H})  $ for \it any \normalfont real valued function $F(\cdot)$.

\par
\begin{figure}
	\centering
	\includegraphics[scale=0.4]{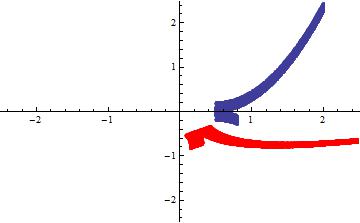}
	\caption{A Rotation}
	\label{sxhmarot}
\end{figure}

\noindent
\noindent
As a third example of transformation,
we act on the cloud of points with a
rotation matrix $C=
\left(\begin{array}{cc}
\cos \displaystyle \frac{\pi}{3} & \sin \displaystyle \frac{\pi}{3} \\
-\sin\displaystyle  \frac{\pi}{3} & \cos \displaystyle \frac{\pi}{3}
\end{array} \right).
$
The initial and the rotated clouds are shown in  Figure \ref{sxhmarot}. The linear coefficients
of the rotated cloud are $\hat{M}=-0.0364518$ and $\hat{H}=0.0143303.$ We found in subsection \ref{sub4},
that a family of cloud invariants, when a cloud is transformed with the rotation matrix $C,$ is $\displaystyle I(M,H)= F \left (\frac{M^2-H}{(H+1)^2} \right ),$ where $F(\cdot)$ is an arbitrary real valued function. Indeed,
we have
$\displaystyle \frac{M^2-H}{(H+1)^2} =\frac{\hat{M}^2-\hat{H}}{(\hat{H}+1)^2} =-0.0126368.$ Henceforth
 we have $I(M,H)= I(\hat{M},\hat{H}),  $ for \it any \normalfont real valued function $F(\cdot)$.
\par

\begin{figure}
\centering
	\includegraphics[scale=0.4]{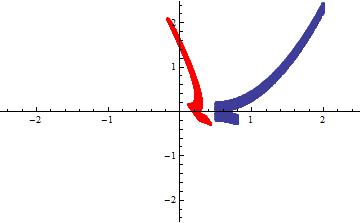}
	\caption{An Arbritrary Matrix}
	\label{sxhmaarb}
\end{figure}

\noindent
Finally, we act on
the cloud of points with
the matrix $D=\left(
\begin{array}{cc}
0.4 & -0.4 \\
-0.05 & 0.9 \\
\end{array}
\right)$. The result
of this transformation
is shown in Figure \ref{sxhmaarb}.
The linear coefficients of the transformed
cloud are $\hat{M}=-4.86159 $ and $\hat{H}=27.4371.$
We found in subsection \ref{linear},
that a family of cloud invariants, when a cloud is transformed with the  matrix $D,$ is $\displaystyle I(M,H)= F \left (\frac{M^2-H}{(8H -1 +10 M)^2} \right ),$ where $F(\cdot)$ is an arbitrary real valued function. Indeed,
we have
$\displaystyle \frac{M^2-H}{(8H -1 +10 M)^2} =\frac{\hat{M}^2-\hat{H}}{(8\hat{H} -1 +10 \hat{M})^2} =-0.000131743.$ Consequently
 we obtain $I(M,H)= I(\hat{M},\hat{H}),  $ for \it any \normalfont real valued function $F(\cdot)$.

\noindent

We note that the results we obtained
by considering the aforementioned cloud
of points verify  our findings in section \ref{invariants}.

\section{Concluding Remarks}
\label{remarks}
We have studied
transformations, with $2 \times 2$
matrices $\left(\begin{array}{cc} \alpha & \beta \\ \gamma & \delta \end{array} \right)$,
$\alpha,\beta,\gamma,\delta \in {\bf R}$, of planar set of points,
called clouds of points for convenience.
Our aim in this paper is to
find cloud invariants, i.e. functions
which take the
same value  when they are evaluated for the
initial and for the transformed cloud of
points.
It is natural the cloud invariants to be
functions of variables which carry information
for the whole cloud.
The cloud invariants we find are functions
of two such variables $M$ and $H.$

\noindent
$M$ and $H$ are functions of the
coordinates of the points of the cloud.
As a result we find that
 any transformation of a cloud
of points by a $2 \times 2$ matrix induces
a nonlinear transformation $(M,H) \rightarrow (\hat{M},\hat{H}),$ $
\hat{M}={\cal M}(M,H,\alpha,\beta,\gamma,\delta),\ \hat{H}={\cal H}(M,H,\alpha,\beta,\gamma,\delta),
$ given explicitly by equations (\ref{transfM}) and (\ref{transfH}), of the variables
$M$ and $H.$

\noindent
$M$ and $H$
originate from
the  best fitting straight line through the cloud of points under consideration. This straight line is  determined by the
least squares fitting technique.
Henceforth by definition a cloud
 invariant is any function $I(M,H)$ which satisfies the relation (\ref{invcon}), $I(M,H)=I(\hat{M},\hat{H}),$  where
 $\hat{M}$ and $\hat{H}$ are the values
 of the variables $M$ and $H$ for the
 transformed cloud.

\noindent
We find cloud invariants by using Lie theory.
Lie theory  replaces the
complicated, nonlinear finite invariance condition (\ref{invcon}) by the  more useful and tractable linear infinitesimal condition (\ref{lin}) provided that the function $I(\hat{M},\hat{H})$ is analytic in the parameters
$\alpha,\beta,\gamma,$ and  $\delta.$ Linear
condition (\ref{lin}) is a set of linear PDEs.
Any solution to this system of PDEs gives a  cloud invariant.

\noindent
Cloud invariants  can be practically
implemented in various fields, e.g. in optical character recognition,  in image analysis and computer graphics techniques, by
providing the necessary tools in order  to identify icons created by the same ``source''.
The cloud invariants we find for the general
four$-$parameter  case, when a cloud is transformed with a matrix $\left(\begin{array}{cc} \alpha & \beta \\ \gamma & \delta \end{array} \right)$, cannot be practically
implemented.

\noindent
However, the cloud invariants we find for various one$-$parameter groups of transformations can be
practically implemented. In particular we find
cloud invariants for a group consisting of diagonal matrices, for a group consisting of upper triangular matrices, for a group consisting of lower triangular matrices, and for the  group of rotations  $SO(2).$

\noindent
More importantly, for the practical implementation of our findings, we find cloud invariants for \it any \normalfont given matrix. We find these cloud invariants by noticing that \it any \normalfont given matrix belongs
to a one$-$parameter ``linear'' set of transformations of the form $\mathcal A + \mathcal B \varphi,$ $\mathcal A$ and  $\mathcal B$ are given $2 \times 2$ matrices, and $\varphi \in {\bf R}.$ Our findings are verified by examples and
simulations in a  cloud of 10.000 points.

\noindent
We expressed the cloud invariants
in
terms of the variables $M$ and $H.$
$M$ and $H$ are essentially the coefficients
of the straight line
which is the best fit for the cloud of points under consideration.
This provides
a natural guide for future research.
With a view to apply our results in fields such as character recognition, the next logical step is to  consider the case where the cloud of points originates from an icon which has a parabolic$-$like shape.

\noindent
In this case we will look for cloud invariants which are  expressed in terms of variables which appear as coefficients,
or variations thereof, of the parabola
which is the best fit for the cloud of points under consideration. For similar reasons, subsequently, we will look for new cloud invariants expressed in terms of variables which appear as coefficients in third or higher degree
curves. We will compare our findings, with
simulations and computational experiments, with
those acquired by other approaches.



\vskip 10 pt

\section{Acknowledgement} The first author would like to express his thanks to Mr. Koutsoulis Nikos, for his attempts to face the problem initially.

\end{document}